\documentclass[11pt]{amsart}
\usepackage{amssymb}
\newtheorem{theorem}{Theorem}[section]
\newtheorem{lemma}[theorem]{Lemma}

\newtheorem{proposition}[theorem]{Proposition}
\newtheorem{corollary}[theorem]{Corollary}
\theoremstyle{definition}

\numberwithin{equation}{section}

\def\IC{{\mathbb C}}
\def\IR{{\mathbb R}}

\def\cN{{\mathcal N}}
\def\cT{{\mathcal T}}

\def\cC{{\mathcal C}}
\def\tr{{\rm tr}\,}
\def\diag{{\rm diag}\,}

\def\Re{{\rm Re}\,}

\def\sige{{\sigma_{\varepsilon}}}
\def\ve{{ \varepsilon }}
\def\Ra{{\Rightarrow}}
\def\[{\left [}
\def\]{\right ]}
\def\({\left (}
\def\){\right )}

\def\Ra{{\ \Rightarrow\ }}

\def\dfrac{\displaystyle\frac}

\def\int{{\bf Int}}

\begin{document}
\openup .9 \jot

\title[ Unitary Similarity Functions
on Lie Products of Matrices]{Preservers of Unitary Similarity Functions \\
on Lie Products of Matrices}

\author{Jianlian Cui}
\address{Department of Mathematics, Tsinghua University, Beijing 100084, P.R. China}
\email{jcui@math.tsinghua.edu.cn}

\author{Chi-Kwong Li}
\address{Department of Mathematics, College of William \& Mary,
Williamsburg, VA 23187, USA}
\email{ckli@math.wm.edu}

\author{Yiu-Tung Poon}
\address{Department of Mathematics, Iowa State University,
Ames, IA 50011, USA}
\email{ytpoon@iastate.edu}

\thanks{}

\dedicatory{}

\commby{}

\begin{abstract}
Denote by $M_n$ the set of $n\times n$ complex matrices.
Let $f: M_n \rightarrow [0,\infty)$ be a continuous map
such that $f(\mu UAU^*)= f(A)$  for any complex unit $\mu$,
$A \in M_n$  and unitary $U \in M_n$, $f(X)=0$ if and only if $X=0$
and the induced map $t \mapsto f(tX)$ is monotonic
increasing on $[0,\infty)$ for any rank 1 nilpotent $X \in M_n$.
Characterizations are  given for surjective maps $\phi$ on $M_n$ satisfying
$f(AB-BA) = f(\phi(A)\phi(B)-\phi(B)\phi(A))$.
The general theorem are then used to deduce results
on special cases when the function is the pseudo spectrum
and the pseudo spectral radius, that answers a question of
Molnar raised at the 2014 CMS summer meeting.
\end{abstract}
\maketitle

{\small
\qquad {\bf AMS Subject Classification} {Primary 15A60, 46B04}

\qquad {\bf Keywords} {Lie product, unitary similarity invariant function, pseudo spectrum.}
}

\section{Introduction}

Let $M_n$ be the set of $n\times n$ matrices.
A function $f: M_n \rightarrow \IR$ is
a {\it radial unitary similarity invariant function} if

\medskip
(P1)  $f(\mu UAU^*)= f(A)$  for a  complex unit $\mu$, $A\in M_n$ and unitary $U\in M_n$.

\medskip\noindent
In \cite{LPS}, the authors studied unitary similarity invariant functions
that are norms on $M_n$, and determine the structure of maps
$\phi: M_n \rightarrow M_n$ satisfying
\begin{equation}\label{eq1}
f(AB-BA) = f(\phi(A)\phi(B) - \phi(B)\phi(A)) \quad \hbox{ for all } A, B \in M_n.
\end{equation}
In \cite[Remark 2.7]{LPS}, it was pointed out that the result actually holds
for more general unitary similarity invariant functions. However, no detail was
given, and it is not   straightforward to apply the results to a specific problem.
For instance, it is unclear how one can apply
the result to study preservers of pseudo spectrum of Lie product of matrices;\footnote{This
is a question  raised by Professor Molnar to the second and third author
at the 2014 Summer Conference of the Canadian Mathematics Society.} see
the definition in Section 3.
To fill this gap, we extend the result in \cite{LPS} to continuous
radial unitary similarity invariant functions  $f: M_n \rightarrow \IR$
satisfying the following properties.

\smallskip
(P2) For any  $X \in M_n$ we have  $f(X) = f(0_n)$ if and only if $X=0_n$, the $n\times n$ zero matrix.

\smallskip
(P3) For any rank 1 nilpotent $X \in M_n$,
the map $t \mapsto f(tX)$ on $[0, \infty)$ is strictly increasing.

\medskip\noindent
For a function $f: M_n \rightarrow [0, \infty)$ satisfying
(P1) -- (P3),  we show that if $\phi: M_n \rightarrow M_n$ is
a surjective map satisfying (\ref{eq1}),  then there is a  unitary $U \in M_n$
and a subset $\cN_n$ of normal matrices in $M_n$ such that $\phi$ has the form
$$\phi(A) = \begin{cases}
\mu_A UA^\dag U^* + \nu_A I_n & A \in M_n \setminus \cN_n \cr
\mu_A U (A^\dag)^* U^* + \nu_A I_n & A \in \cN_n,
\end{cases}$$
where

(a) $\mu_A, \nu_A\in \IC$ with $|\mu_A| = 1$, depending  on $A$,

(b) $A^\dag = A$, $\overline A$, $A^t$ or $A^*$, and

(c) $\cN_n$ depending on the given unitarily invariant function $f$.

\medskip\noindent
The proof of this result will be given in Section 2.
In Section 3, we apply the main result to the case when $f$ is the
pseudo spectral radius, and then obtain the result for the case when
$f$ is the pseudo spectrum.

\medskip
For other preserver problems on different types of products on matrices and operators,
one may see \cite{CLS,CLS2,CFL,GL,LT} and their references.

\section{Main theorem}

In this section, we  prove Theorem \ref{main0} extending
the result in \cite{LPS}.
We use similar ideas in \cite{LPS} with some intricate
arguments to make the extension possible.

\begin{theorem} \label{main0}
Let $f: M_n \rightarrow [0,\infty)$ be a function on $M_n$
satisfying {\rm (P1) -- (P3)}.
Suppose $n \ge 3$, and $\phi: M_n \rightarrow M_n$ is a surjective map
satisfying
$$f([\phi(A),\phi(B)]) = f([A,B]).$$
Then there is a unitary matrix $U$ and a subset $\cN_n$ of normal matrices with non-collinear eigenvalues
such that $\phi$ has the form
$$\phi(A) = \begin{cases}
\mu_A U\tau(A) U^* + \nu_A I_n & A \in M_n \setminus \cN_n \cr
\mu_A U \tau(A)^* U^* + \nu_A I_n & A \in \cN_n,
\end{cases}$$
where $\mu_A, \nu_A\in \IC$ with $|\mu_A| = 1$ depending on $A$,
and $\tau$ is one of the maps: $A \mapsto A$, $A \mapsto \overline A$,
$A \mapsto A^t$ or $A \mapsto A^*$.
\end{theorem}

To prove the above theorem, we need the following result from \v
Semrl \cite{Sem}.

\begin{theorem} \label{semrl}
Suppose $n \ge 3$, and $\phi: M_n \rightarrow M_n$ is a bijective map
satisfying $$[A,B] = 0_n \quad \iff \quad [\phi(A),\phi(B)] = 0_n.$$
Let $\Gamma$ be the set of matrices $A$ such that the Jordan form of
$A$ only has Jordan blocks of sizes 1 or 2. Then there are an
invertible matrix $S$, an automorphism $\tau$ of the complex field
and a regular locally polynomial map $A \mapsto p_A(A)$ such that
\begin{eqnarray}\label{semrl.1}
\phi(A) = S(p_A(A_\tau^\dag))S^{-1} \quad\hbox{for all } A\in \Gamma.
\end{eqnarray}
Here, $X_{\tau}$ is the matrix whose $(i,j)$-entry is
$\tau(X_{ij})$, and $A^\dag = A$ or $A^t$.
\end{theorem}

Our proof strategy is to show that $\phi(A)$ has the asserted
form described in the theorem for a special class $\cC_1$ of
matrices $A$. Then we modify the map $\phi$ to
$\phi_1$ so that it will satisfy the same hypothesis of $\phi$ with
the additional assumption that $\phi(X) = X$ for every $X \in \cC_1$.
Then we can set $B = \phi(A)$ for a certain matrix $A$ not in $\cC_1$
and use the condition that
$$f([A,X]) = f([\phi_1(A),\phi_1(X)]) = f([B,X]) \quad \hbox{ for all } X \in \cC_1$$
to show that $B = \phi_1(A)$ also has the asserted form. Thus,
$\phi_1$ has the asserted form for a larger class $\cC_2$
of matrices, and so on and so forth until we show that the modified
map will fix every matrices after a finite number of steps.

In the next few lemmas, we will focus on the conditions the relations between
a pair of matrices $A$ and $B$ such that
$$f([A,X]) = f([B,X]) \quad \hbox{ for all } X \in \cC$$
for a certain subset $\cC$ of matrices.

In using the conditions (P1) and (P3), we note that every rank one nilpotent matrix is of the
form $xy^*$ for some non-zero orthogonal vectors, $x$ and $y$. Also, $xy^*$ is unitarily similar
to $\|x\|\|y\|E_{21}$.
These facts will be used frequently in our proofs.

Denote by $\sigma(A)$ the spectrum of $A$ and $N(A)$ the null space of $A$.

\begin{lemma}\label{lem.1}
For any two matrices $A$ and $B$, if
\begin{eqnarray}\label{AB}
f([A,X]) = f([B,X]) \quad\hbox{for all rank one }X\in M_n,
\end{eqnarray}
then there are $\mu, \nu\in\IC$ with $|\mu| = 1$ such that
one of the following holds with $\hat A = \mu A + \nu I_n$.
\begin{enumerate}
\item[\rm (a)] $\sigma(B) = \sigma(\hat A)$ and
for any $\lambda\in \sigma(\hat A)$,
$$N(B -\lambda I_n) = N(\hat A-\lambda I_n)
\quad\hbox{and}\quad N(B^t -\lambda I_n) = N(\hat A^t-\lambda I_n).$$

\item[\rm (b)] The eigenvalues of $A$ are not collinear,
$\sigma(B) = \overline {\sigma(\hat A)}$
and for any $\lambda\in \sigma(\hat A)$,
$$N(B -\overline \lambda I_n) = N(\hat A-\lambda I_n)
\quad\hbox{and}\quad
N(B^t -\overline \lambda I_n) = N(\hat A^t-\lambda I_n).$$
\end{enumerate}
\end{lemma}

\it Proof. \rm Note that for any rank one matrix $X = xy^t$, $[C,X] = 0$
if and only if $x$ and $y^t$ are the right and left
eigenvectors of $C$ corresponding to the same eigenvalue. To
see this, as $[C,X]  = (Cx)y^t - x(y^tC)$, then $[C,X] = 0$ if
and only if $Cx = \lambda x$ and $y^tC = \lambda y^t$ for some
$\lambda\in \IC$.

Suppose $A$ and $B$ satisfy (\ref{AB}). By the above observation on
rank one matrices and property (P2) of $f$,
$A$ and $B$ must have the same set of left and
right eigenvectors. Furthermore, $x_1$ and $x_2$ are the right
eigenvectors of $A$ corresponding to the same eigenvalue if and only
if the two eigenvectors correspond to the same eigenvalue of $B$.
Thus, the eigenvalues of $A$ and $B$ have the same geometric
multiplicity.

Let $\lambda_1,\dots,\lambda_k$ be the distinct eigenvalues of $A$
with $x_1,\dots,x_k$ and $y_1,\dots, y_k$ being the right and left
eigenvectors. Also for each pair of eigenvectors $x_i$ and $y_i^t$,
let $\gamma_i$ be the corresponding eigenvalue of $B$.
Take $X_{ij} = x_iy_j^t$. Then $AX_{ij} = \lambda_i X_{ij}$ and $X_{ij} A =
\lambda_j X_{ij}$. Using (P1), we see that for any $1\le i,j\le n$,
$$ f([A,X_{ij}]) = f(\lambda_i X_{ij} - \lambda_j X_{ij})
= f((\lambda_i-\lambda_j) X_{ij}) = f(|\lambda_i-\lambda_j| X_{ij}).
$$
Similarly,
$f([B,X_{ij}]) = f((\gamma_i-\gamma_j)X_{ij})
= f(|\gamma_i-\gamma_j| X_{ij}).$

By the fact that $f([A,X_{ij}]) = f([B,X_{ij}])$ and Property (P3),
$$|\lambda_i-\lambda_j| = |\gamma_i-\gamma_j| \quad \hbox{ for all }
1\le i,j\le k.$$
As a result, there are
$\mu,\nu\in \IC$ with $|\mu|=1$ such that either
\begin{enumerate}
\item $\gamma_i = \mu \lambda_i + \nu$ for all $1\le i\le k$; or
\item the eigenvalues of $A$ are non-collinear and
$\overline \gamma_i = \mu \lambda_i + \nu$ for all $1\le i\le k$.
\end{enumerate}
Then the result follows with $\hat A = \mu A + \nu I_n$.
\qed

\begin{lemma}\label{lem.2}
Suppose $A$ and $B$ commute and satisfy {\rm(\ref{AB})}.
If $A$ has at least two distinct eigenvalues, then
there are $\mu,\nu \in \IC$ with $|\mu| = 1$ such that
\begin{enumerate}
\item[\rm (a)] $B = \mu A + \nu I_n$, or
\item[\rm (b)] $A$ is normal with non-collinear eigenvalues
and $B = \mu A^* + \nu I_n$.
\end{enumerate}
\end{lemma}

\it Proof. \rm
As $A$ and $B$ commute, there is a unitary matrix $U$
such that both $U^*AU$ and $U^*BU$ are upper triangular, see \cite[Theorem 2.3.3]{HJ}.
Replacing  $(A,B)$ with $(U^*AU,U^*BU)$,
we may assume that $A$ and $B$ are upper triangular.

As $A$ and $B$ satisfy (\ref{AB}), Lemma \ref{lem.1} holds.
Suppose Lemma \ref{lem.1}(a) holds with $\hat A = \mu A + \nu I_n$.
Notice that $\sigma(B) = \sigma(\hat A)$ and
$$f([\hat A,X]) = f([\mu A + \nu I_n,X ]) =
f([B,X]) \quad\hbox{for all rank one } X\in M_n.$$
Suppose $\lambda$ is an eigenvalue of $\hat A$ and $y \in N(\hat A^t
- \lambda I_n)$. For any $z\in \IC^n$, let $Z = z y^t$. Then $Z \hat
A = \lambda Z$ and $[\hat A,Z] = (\hat A - \lambda I_n)Z$. Note that
$(\hat A - \lambda I_n)Z$ has rank at most one and $\tr((\hat A-
\lambda I_n)Z) = \tr ([\hat A,Z]) = 0$, so $(\hat A- \lambda I_n) Z$
is unitarily similar to $\|(\hat A - \lambda I_n)z\|\|y^t\| E_{12}$.
Thus,
$$f([\hat A,Z]) = f(\|(\hat A - \lambda I_n)z\|\, \|y^t\| E_{12}).$$
Similarly, $f([B,Z]) =
f(\|(B - \lambda I_n)z\|\, \|y^t\| E_{12})$. Hence, by (P1) and (P3),
$$\|(\hat A - \lambda I_n)z\| = \|(B - \lambda I_n)z\|
\quad\hbox{for all }z \in \IC^n \hbox{ and } \lambda\in \sigma(\hat A).$$
As a result,
\begin{eqnarray*}
&& z^*\hat A^*\hat A z - 2\Re (\overline \lambda z^*\hat Az)
+ |\lambda|^2z^* z = \|(\hat A - \lambda I_n) z\|^2 \\
&=& \|(B - \lambda I_n) z\|^2
= z^*B^*B z - 2\Re (\overline \lambda z^*Bz) + |\lambda|^2z^* z.
\end{eqnarray*}
This implies that
$$2\Re( \overline\lambda z^*(\hat A-B) z ) = z^*(\hat A^*\hat A - B^*B)z
\quad\hbox{for all }z \in \IC^n \hbox{ and } \lambda\in \sigma(\hat
A).$$ As $A$ has at least two distinct eigenvalues, so does $\hat
A$. Taking any $\lambda,\gamma\in \sigma(\hat A)$ with $\lambda\ne
\gamma$, we have
$$2\Re( \overline\lambda z^*(\hat A-B) z ) = z^*(\hat A^*\hat A - B^*B)z
= 2\Re( \overline\gamma z^*(\hat A-B) z ).$$
Thus, $W((\overline{\lambda-\gamma}) (\hat A-B)) \subseteq i\IR$, where
$W(X)$ is the numerical range of $X$. Then $(\overline{\lambda -
\gamma}) (\hat A-B)$ is a skew-Hermitian matrix and hence $\hat A-B$
is a diagonal matrix. Now for any $1\le i \le n$, $b_{ii} \in
\sigma(B) = \sigma(\hat A)$. Also the $i$th entry of $(B - b_{ii}
I_n) e_i$ is zero while only the $i$th entry of $(\hat A - B) e_i$
can be nonzero. Then
\begin{eqnarray*}
\|(B- b_{ii} I_n) e_i\|^2
&=& \|(\hat A -  b_{ii} I_n) e_i \|^2
= \|(B- b_{ii} I_n) e_i + (\hat A - B) e_i\|^2 \cr
&=& \|(B- b_{ii} I_n) e_i\|^2 + \|(\hat A - B) e_i\|^2.
\end{eqnarray*}
Thus, $(\hat A - B)e_i = 0$ for all $1\le i\le n$ and hence $B = \hat A$.

Now suppose Lemma \ref{lem.1}(b) holds. Then by a
similar argument, we can show that
\begin{eqnarray}\label{eq.101}
\|(\hat A - \lambda I_n)z\| = \|(B - \overline \lambda I_n)z\|
\quad\hbox{for all }\lambda\in \sigma(\hat A) \hbox{ and } z \in \IC^n
\end{eqnarray}
and so $(\overline{\lambda - \gamma})\hat A - ( \lambda- \gamma ) B$
is a skew-Hermitian matrix. It follows that
$(\overline{\lambda - \gamma}) T_A
- ( \lambda- \gamma ) T_B = 0$, or equivalently, $T_B =
\frac{\overline{\lambda - \gamma}}{\lambda - \gamma} T_A$, where
$T_A$ and $T_B$ are the strictly upper triangular parts of $A$ and
$B$. Now as the eigenvalues of $A$ and hence $\hat A$ are not
collinear, we can always find another $\omega\in \sigma(\hat A)$
such that $\frac{\overline{\lambda - \omega}}{\lambda - \omega} \ne
\frac{\overline{\lambda - \gamma}}{\lambda - \gamma}$. Then the
above equation is possible only if $T_A = T_B = 0$. In this case,
$A$ and $B$ are both diagonal and hence normal. Then (\ref{eq.101})
implies that $\hat A = \overline B$. \qed

From Lemma \ref{lem.2}, we have the following consequence for
diagonalizable matrices.

\begin{corollary}\label{cor.1}
Suppose $A$ and $B$ satisfy (\ref{AB}) and $A$ is diagonalizable.
Then there are $\mu,\nu\in\IC$ with $|\mu| = 1$ such that
\begin{enumerate}
\item[\rm (a)] $B = \mu A + \nu I_n$, or
\item[\rm (b)] $A$ is normal with non-collinear eigenvalues and $B = \mu A^* + \nu I_n$.
\end{enumerate}
\end{corollary}

\it Proof. \rm Suppose $A$ is diagonalizable.
Then $A = SDS^{-1}$ for some invertible $S$ and diagonal $D$.
By Lemma \ref{lem.1}, $B = S(\mu D + \nu I_n)S^{-1}$
or $B = S(\mu \overline D + \nu I_n)S^{-1}$.
If $A$ has only one eigenvalue, then $A$ is a scalar matrix
and so is $B$. Then the result follows.
Suppose $A$ has at least two eigenvalues.
As $A$ and $B$ commute, the result now follows by Lemma \ref{lem.2}.
\qed

\begin{lemma}\label{lem.3}
For any two matrices $A$ and $B$, if
\begin{eqnarray}\label{AB2}
f([A,X]) = f([B,X])\quad\hbox{for all }X\in M_n,
\end{eqnarray}
then there are $\mu, \nu\in\IC$ with $|\mu| = 1$ such that
\begin{enumerate}
\item[\rm (a)] $B = \mu A + \nu I_n$, or
\item[\rm (b)] $A$ is normal with non-collinear eigenvalues and $B = \mu A^* + \nu I_n$.
\end{enumerate}
\end{lemma}

\it Proof. \rm
Suppose $A$ and $B$ satisfy (\ref{AB2}).
Then clearly $A$ and $B$ commute.
If $A$ has at least two eigenvalues, then the result
follows from Lemma \ref{lem.2}.

Suppose $A$ has only one eigenvalue, say $\lambda$. Then by Lemma
\ref{lem.1}, $B$ has one eigenvalue only, say $\gamma$. Write $A =
SJS^{-1} + \lambda I_n$, where $S$ is invertible and $J = J_{n_1}
\oplus \cdots \oplus J_{n_s}$ is the Jordan form of $A$ with $n_1\ge
\cdots \ge n_s$. Now as $A$ and $B$ satisfy (\ref{AB2}), $A$ and $B$
have the same set of commuting matrices. Then $B = Sp(J)S^{-1} +
\gamma I_n$ for some polynomial $p$ of degree at most $m=n_1 - 1$
with $p(0) = 0$.

By a similar argument as in Lemma \ref{lem.2}, we can show that
$$\|(B - \gamma I_n)z\| = \|(A - \lambda I_n)z\| \quad \hbox{for all }z\in \IC^n.$$
Then there is a unitary matrix $W$ such that
$$Sp(J)S^{-1} = (B - \gamma I_n) = W(A - \lambda I_n) = WSJS^{-1}.$$
Write $S = UT$ for unitary $U$ and upper triangular $T$,
$V = U^*WU$ and $p(x) = \sum_{i=1}^m c_i x^i$. Then we have
\begin{eqnarray}\label{eq.1}
Tp(J)T^{-1} = VTJT^{-1}.
\end{eqnarray}
Notice that both $Tp(J)T^{-1}$ and $TJT^{-1}$ are strictly upper
triangular. Furthermore, the first $n_1 - 1$ entries in the
super-diagonal of $Tp(J)T^{-1}$ are $c_1$ times the corresponding
$n_1 - 1$ super-diagonal entries of $TJT^{-1}$.

As $V$ is unitary, we must have $|c_1| = 1$ and $V = c_1 I_{n_1 - 1}
\oplus V_1$ for some unitary $V_1 \in M_{n-n_1+1}$. Now comparing
the first $n_1\times n_1$ principal submatrices in (\ref{eq.1}), we
have
$$T_1 p(J_{n_1}) T_1^{-1}
= (c_1 I_{n_1-1} \oplus [v_{n_1,n_1}] )T_1 J_{n_1}T_1^{-1} = c_1 T_1
J_{n_1}T_1^{-1},$$ where $T_1$ is the $n_1\times n_1$ principal
submatrix of $T$. Therefore, $T_1 \left(\sum_{i=2}^m c_i J_{n_1}^i
\right) T_1^{-1}= 0$ and so $\sum_{i=2}^m c_i J_{n_1}^i = 0$. Hence,
$c_2 = \cdots = c_m = 0$. Then $p(x) = c_1 x$ and so $B = c_1 A +
(\gamma - c_1\lambda)I_n$. \qed

\medskip
We are now ready to present the following.

\medskip\noindent
\it Proof of Theorem \ref{main0}. \rm

{\bf First we assume that $\phi$ is bijective.}
Suppose $\phi$ is a bijective map satisfying
$$f([A,B]) = f([\phi(A),\phi(B)]) \quad\hbox{for all }A,B\in M_n.$$
Because $f(X) = f(0)$ if and only if $X = 0$ by (P2),
we see that $[A,B] = 0$ if and only if $[\phi(A),\phi(B)] = 0$.
We can apply  Theorem \ref{semrl} and conclude that
$\phi$ has the form (\ref{semrl.1})
with $A^\dag = A$ or $A^t$.
In particular, for any rank one matrix $R \in M_n$,
there are $\mu_R, \nu_R\in \IC$ such that
$$\phi(R) = S(\mu_R R_\tau^\dag + \nu_R I_n) S^{-1}.$$
Without loss of generality, we may assume that
$\mu_R > 0$ and $\nu_R = 0$.

Here we consider only the case when $A^\dag = A$.
The case when $A^\dag = A^t$ is similar.
Fixed an orthonormal basis $\{x_1,\dots,x_n\}$ and define $X_{ij} = x_i x_j^*$.
Take $\alpha = (\alpha_1,\dots,\alpha_n)\in\IC^n$
and let $A = \sum_{j = 1}^n \alpha_j X_{j1}$.
For $k = 2,\dots,n$,
\begin{equation} \label{2001}
 f(\mu_A \mu_{X_{kk}} \tau(\alpha_k) S(X_{k1})_\tau S^{-1})
= f([\phi(A),\phi(X_{kk})])
= f([A,X_{kk}])
= f(\alpha_k X_{k1}).
\end{equation}
In particular,
if $Z = \mu_A \mu_{X_{22}} S(X_{21})_\tau S^{-1}$, then
$$f(\tau(\alpha)  Z) = f(\alpha X_{21}) \quad \hbox{ for all } \alpha\in \IC.$$
Since the induced maps $g: t \mapsto f(tX_{21})$ and
$h: t \mapsto f(t Z)$ are monotonic
increasing for $t \in [0, \infty)$ by (P3), we have a well defined continuous map
$\tau(\alpha) = h^{-1} g(\alpha)$ on an interval $[0, d)$ for some  $d > 0$.
Since $\tau$ is an automorphism on $\IC$,
it is either the identity map $\lambda \mapsto \lambda$ or
the conjugate map $\lambda \mapsto \overline \lambda$; for example, see \cite{Kest}.

Furthermore, as $f([X_{32},X_{22}]) = f(X_{32}) = f([X_{32},X_{33}])$,
\begin{eqnarray*}
&&f(\mu_{X_{32}} \mu_{X_{22}} S(X_{32})_\tau S^{-1}) =
f([\phi(X_{32}),\phi(X_{22})]) \cr
&=& f([\phi(X_{32}),\phi(X_{33})])
= f(\mu_{X_{32}} \mu_{X_{33}} S(X_{32})_\tau S^{-1}).
\end{eqnarray*}
Thus, $\mu_{X_{22}} = \mu_{X_{33}}$ by (P3).
By (\ref{2001}) and the fact that $f(\xi X_{21}) = f(\xi X_{31})$ for
all $\xi \in \IC$, we have
$$f(S(X_{21})_\tau S^{-1}) = f(S(X_{31})_\tau S^{-1}).$$

We now claim that $S$ is a multiple of some unitary matrix.
If not, then there is a pair of orthonormal vectors $y_2,y_3$
such that $\|Sy_2\| \ne \|Sy_3\|$.
Extend $y_2,y_3$ to an orthonormal basis $\{y_1,y_2,y_3,\dots,y_n\}$
and let $x_j = (y_j)_{\tau^{-1}}$. Then $\{x_1,\dots,x_n\}$ also forms
an orthonormal basis. By the above study, we have
$$f(\|Sy_2\| \|y_1^* S^{-1}\|\, E_{12})
= f(S(X_{21})_\tau S^{-1})
= f(S(X_{31})_\tau S^{-1})
= f(\|Sy_3\| \|y_1^* S^{-1}\| \, E_{12}),$$
which contradicts that $\|Sy_2\| \ne \|Sy_3\|$.
Thus, $S$ is a multiple of some unitary matrix.
By absorbing the constant term,
we may assume that $S$ is unitary.
Now for any rank one matrices $R$ and $S$,
$$f([R,S]) = f([\phi(R),\phi(S)])
= f(\mu_R\mu_S [R_\tau,S_\tau]).$$
By (P1), $f([R,S]) = f([R_\tau, S_\tau])$ whenever $[R,S]$
is a rank one nilpotent matrix, and hence $\mu_R\mu_S = 1$ in this
case.

Now for any rank one matrix $A$, we can always find two other rank
one matrices $B$ and $C$ such that $[A,B]$, $[A,C]$ and $[B,C]$ are
all rank one nilpotents. Then we must have $\mu_A\mu_B = \mu_A\mu_C
= \mu_B\mu_C = 1$. As all $\mu_A,\mu_B,\mu_C$ are positive real
numbers, the equality is possible only when $\mu_A = \mu_B = \mu_C =
1$. Then we have $\phi(A) = SA_\tau S^{-1} = SA_\tau S^*$ for
all rank one $A$.

By replacing $\phi$ with the map $A \mapsto S^*\phi(A) S$, we may
assume that $\phi(X) = X^+$ for all rank one matrices $X$, where
$X^+ = X$, $\overline X$, $X^t$ or $X^*$. Then
$$f([A,B]) = f([\phi(A),\phi(B)])
= f([A^+, B^+]) = f([A, B]^+)$$
for all rank one $A,B\in M_n.$
Notice that the set
$$\{X: X = [A,B] \hbox{ for some rank one $A$ and $B$} \}$$
contains the set of trace zero non-nilpotent
matrices with rank at most two and so is dense in the set of trace
zero matrices with rank at most two. Thus, by continuity of $f$ we see that
$$f(X) = f(X^+)\quad  \hbox{for all trace zero matrices $X$ with rank at most two}.$$
Now define $\Phi:M_n \to M_n$ by $A\mapsto \phi(A)^+$. Then $\Phi(X)
= X$ for all rank one matrices $X$. For any $A\in M_n$ and rank one
matrix $X\in M_n$, as $[A,X]$ is a trace zero matrix with rank at
most two,
$$f([A,X]) = f([\phi(A),\phi(X)])
= f([\phi(A),X^+]) = f([\phi(A)^+,X]) = f([\Phi(A),X]).$$
Thus,
$f([A,X]) = f([\Phi(A),X])$ for all rank one $X$. Then Corollary
\ref{cor.1} implies that $\Phi(A) = \mu_A A + \nu_A I_n$ or $\Phi(A)
= \mu_A A^* + \nu_A I_n$ for all diagonalizable matrices $A$ and the
latter case happens only when $A$ is normal with non-collinear eigenvalues.

After absorbing the constants $\mu_A$ and $\nu_A$, we may assume
that $\Phi(X) = X$ for all non-normal diagonalizable matrices $X$.
Then
$$f([A,B]) = f([\phi(A),\phi(B)])= f([\Phi(A),\Phi(B)]^+) = f([A,B]^+)$$
for all non-normal diagonalizable matrices $A$ and $B$. Since the
set of all non-normal diagonalizable matrices is dense in $M_n$, we
see that $f([A,B]) = f([A,B]^+)$ for all $A,B\in M_n$. Then for
any $A\in M_n$,
$$f([A,X]) = f([\phi(A),\phi(X)]) = f([\Phi(A),\Phi(X)]^+) = f([\Phi(A),X])$$
for all non-normal diagonalizable matrices $X$,
and so $f([A,X]) = f([\Phi(A),X])$ for all $X\in M_n$ by the continuity of $f$. Now the
result follows by Lemma \ref{lem.3}.

\medskip
{\bf Finally, we show that one only needs the surjective
assumption on $\phi$}.
For any $A, B\in M_n$, we say $A \sim B$ if
$$f([A,X]) = f([B,X])\quad\hbox{for all }X\in M_n.$$
Clearly, $\sim$ is an equivalence relation and for each $A\in M_n$,
denote by $S_A = \{B: B\sim A\}$ the equivalence class of $A$. By
Lemma \ref{lem.3}, either
\begin{enumerate}
\item[\rm (I)] $S_A$ is the set of matrices of the form
$\mu A + \nu I$ for some $\mu,\nu\in \IC$ with $|\mu| = 1$, or
\item[\rm (II)] $A$ is normal and $A \sim A^*$,
$S_A$ is the set of matrices of the form
$\mu A + \nu I$ or $\mu A^*+ \nu I$ for some $\mu,\nu\in \IC$ with $|\mu| = 1$.
\end{enumerate}

Pick a representative for each equivalence class and write
$\mathcal{A}$ for the set of these representatives.  Since $\phi$ is
surjective, $S_A$ and $\phi^{-1}(S_A)$ have the same cardinality $c$
for every $A \in \mathcal{A}$.  Thus there exists a map $\psi:M_n
\to M_n$ which maps $\phi^{-1}(S_A)$ bijectively onto $S_A$ for each
$A \in \mathcal{A}$.  Clearly $\psi$ is bijective and $\psi(A) \sim
\phi(A)$ for all $A \in M_n$. Then, for any $A,B\in M_n$,
$$f([A,B]) = f([\phi(A),\phi(B)]) = f([\psi(A),\phi(B)]) =
f([\psi(A),\psi(B)]).$$
That is, $\psi$ is bijective map satisfying
(\ref{AB}). By Theorem \ref{main0}, $\psi$ has the desired form and
hence so does $\phi$, as $\psi(A) \sim \phi(A)$ implies $\phi(A) =
\mu \psi(A) + \nu I$ or $\phi(A) = \mu \psi(A)^* + \nu I$ when
$\psi(A)^*$ is normal and $\psi(A)^* \sim \psi(A)$. \qed

\medskip\noindent
{\bf Remark} Using the argument in the last part of the proof on the replacement of the bijective assumption by the surjective assumption on $\phi$, one may further weaken the surjective assumption on $\phi$ by any one of
the following (weaker) assumptions on the following modified map
$\tilde \phi$ defined by
$$\tilde \phi(X) = \phi(X) - \tr(\phi(X))I/n$$
on the set $M_n^0$ of trace zero matrices in $M_n$.

(a) The map $\tilde \phi: M_n^0 \rightarrow M_n^0$ surjective.

(b) For any $A \in M_n^0$ the range of $\tilde \phi$ contains a matrix of the form
$e^{it} A$ for some $t \in [0, 2\pi)$.

\section{Pseudo spectrum and pseudo spectral radius}

In this section, we use Theorem \ref{main0} to study maps preserving the
pseudo spectral radius (see the definitions below) of the Lie product of matrices.
Then we further deduce
the result for maps preserving the pseudo spectrum. As one shall see,
with considerable effort,
one will be able to get more specific structure of the preserving maps.

For $\varepsilon >0$, define the $\varepsilon$-pseudospectrum
$\sigma _\varepsilon (A)$ of $A\in M_n$ as
$$\sigma _\varepsilon (A)
= \{z \in \sigma(A+E): E \in M_n, \ \|E\| < \varepsilon\}
= \{z\in \IC: s_n(A - zI_n) < \varepsilon\},
$$
where $s_1(X) \ge \cdots \ge s_n(X)$ denote the singular values of
$X \in M_n$,
and the $\varepsilon$-pseudospectral radius
$r _\varepsilon (A)$ of $A\in M_n$ as
$$r _\varepsilon (A) = \sup\{ |\mu|: \mu \in \sigma _\varepsilon (A)\}.$$
Note that the pseudo spectral radius is useful in studying the stability of
matrices under perturbations, and there are efficient algorithm for its
computation; see for example, \cite{O1} and its references.
Preservers of pseudo spectrum has been considered for several types of products in \cite{Cui1} (see also \cite{Cui2}).
Here we  characterize the preservers of pseudo spectral radius and pseudo
spectrum for Lie products.  We first prove the following.

\begin{theorem} \label{3.1} Suppose $n \ge 3$ and $\varepsilon >0$.
Then a surjective map $\phi: M_n \rightarrow M_n$ satisfies
$$r _\varepsilon([A,B])=r _\varepsilon ([\phi (A),\phi (B)])
\qquad\hbox{ for all } A, B \in M_n$$
if and only if there is a unitary $U \in M_n$ such that
$$\phi(A) =
\mu_A U\tau(A) U^* + \nu_A I_n  \quad \hbox{ for all } A \in M_n,$$
where $\mu_A, \nu_A\in \IC$ with $|\mu_A| = 1$, depending on $A$,
and $\tau$ is one of the  following maps: $A\mapsto A$, $A \mapsto \overline A$,
$A \mapsto A^t$ or $A \mapsto A^*$.
\end{theorem}

\it Proof. \rm The sufficiency can be readily checked. To prove the necessity, Let $f(A)=r _\varepsilon (A)$ for $A\in M_n$.
It is clear  that $f$ is a continuous map satisfying (P1) and (P2). Suppose $X$ is a rank one nilpotent matrix. It follows from Proposition 2.4 in  \cite{Cui2} that  $r_{\ve}(X)= \sqrt{\ve^2+\|X\|\ve} $. Hence,  (P3) is also satisfied.
So, we can apply Theorem \ref{main0} and conclude that $\phi$ has the form
in Theorem \ref{main0}.
To get the desired conclusion, we need to show that the set $\cN$
is empty. Assume not, and there is $A \in \cN$.
Since $A$ is normal with non-collinear eigenvalues, there is a unitary $V$ and  $\gamma, \xi \in \IC$ such that
$$V(\tau(A)-\xi I) V^* = \gamma \diag(1, \mu, 0, \mu_4, \dots, \mu_n),$$
where $\mu \notin \IR$.
Let $B \in M_n$ be such that
$$\tilde B =
V\tau(B)V^* = \left[\begin{matrix} 0 & 1 & 0 \cr a & 0 & b \cr 0 & c & 0 \cr
\end{matrix}\right] \oplus O_{n-3},$$
where $a = (1-\bar\mu)/(1-\mu)$, $b > 0$  and $c = b \bar \mu/\mu$.
Then
$$\tilde B \tilde B^* = \left[\begin{matrix}
1 & 0 & \bar c\cr 0 & |a|^2 + |b|^2 & 0 \cr c & 0 & |c|^2 \cr
\end{matrix}\right] \quad \hbox{ and }
\quad
\tilde B^* \tilde B = \left[\begin{matrix}
|a|^2 & 0 & \bar a b\cr 0 & 1 + |c|^2 & 0 \cr \bar b a  & 0 & |b|^2
\cr\end{matrix}\right]$$
and we can choose $b > 0$ so that $\tilde B$ is not normal, and neither is $B$.
As a result, $\phi(B) = \mu_B U\tau(B)U^* + \nu_B I$.

Now,
$$C_1 = V[\tau(A),\tau(B)]V^*
= \gamma \left[\begin{matrix} 0 & 1-\mu & 0 \cr \bar \mu-1  & 0 & b\mu \cr
0 & -b \bar \mu & 0 \cr
\end{matrix}\right] \oplus O_{n-3}
$$
is normal with eigenvalues $s_\pm = \pm \gamma \sqrt{|1-\mu|^2 + b^2|\mu|^2}$
so that
$$r_\varepsilon([A,B]) = r_\varepsilon([\tau(A), \tau(B)]) =
|\gamma|\sqrt{|1-\mu|^2 + b^2|\mu|^2} + \varepsilon.$$
However, $[\phi(A),\phi(B)]$ is unitarily similar to
$$C_2 =\mu _A\mu _B \bar\gamma
\left[\begin{matrix} 0 & 1-\mu & 0 \cr
(1-\bar{\mu})^2/(\mu -1) & 0 & b\bar\mu \cr 0 &  -b\bar{\mu}^2/\mu & 0 \cr
\end{matrix}\right] \oplus O_{n-3}.
$$
One readily checks that the matrix $C_2$ is normal if and only if
$\mu $  is pure imaginary.
In all other cases,
there is a unitary $R \in M_n$ obtained from
$I_n$ by changing the $(1,1), (1,3), (3,1), (3,3)$ entries so that
$$RC_2R^* = \bar \gamma
\left[\begin{matrix} 0 & c_1 & 0 \cr
c_2 & 0 & c_3 \cr 0 & 0 & 0 \cr
\end{matrix}\right] \oplus O_{n-3}.$$
If $C_2$ has singular values $s_1  \ge s_2$, then
$$|\gamma |^2(|c_1|^2 + |c_2|^2 + |c_3|^2)
= \tr(C_2 C_2^*) = \tr(C_1C_1^*) = |\gamma |^2(s_+^2 + s_-^2).$$
Because $C_2$ is not normal, $s_1 < s_+$, we see that $s_2 > s_-$.
Then for any $z \in \IC$, if
$\tilde C - z I$ has singular values $s_1(z)\ge s_2(z)$, then
$$s_1(z)^2 + s_2(z)^2  = 2|z|^2 + |c_1|^2 + |c_2|^2 + |c_3|^2
2|z|^2 + s_+^2 + s_-^2 = s_+(z)^2 + s_-(z)^2,$$
where $s_+(z) \ge s_-(z)$ are the singular values of
$C_1 - zI$.
Again, because $C_2-zI$ is not normal, we see that
$s_+(z) > s_1(z) \ge s_2(z) > s_-(z)$.
It follows that $s_2(z) > s_-(z)$
for any $z \in \IC$ with
$|z| \le |\gamma|\sqrt{|1+\mu|^2 + b^2|\mu|^2} + \varepsilon]$.
Thus,
$$\max\{ z \in \IC: s_2(C_2 - zI) \le \varepsilon\}
< \max\{ z \in \IC: s_2(C_1 - zI) \le \varepsilon\}.$$
So, if a normal matrix $A$ with three collinear eigenvalues
$\gamma + \nu, \gamma\mu  + \nu, \nu$ so that  $\mu$ is not real
and $\mu \ne \pm i$, then $A \notin \cN$.
Clearly, if $A \in \cN$  has eigenvalues of the form
$\gamma  + \nu, \gamma +i \nu, \gamma$,
then $\tau(A)^*$ can be viewed as a multiple of
$\tau(A)$. Thus, we may assume that $A\notin \cN$ by adjusting
$\mu_A$ and $\nu_A$.
The result follows.
\qed

We will use the above theorem to determine the structure of preservers
of the pseudo spectrum of Lie product of matrices.
To achieve this, we need a characterization of
normal matrices $A$ with two distinct eigenvalues, i.e.,
$A - bI$ is a nonzero multiple of a
rank $k$ orthogonal projection $P$  with $1 \le k < n$; see Proposition
\ref{3.2a} below. The proof depends on the following lemma.

\begin{lemma} \label{lemt} Suppose $C =C_1\oplus O_{n-3}$,
where $C_1\in M_3$ has rank $\le 2$ and $\tr C_1=0$. Then
for every $\varepsilon >0$, $\sige(C)=\sige(C_1)$. Furthermore, suppose for
$t\in \IR$,  $$f(\lambda,t)=\det(\lambda I_3-(C_1-tI_3)^*(C_1-tI_3))
=\lambda^3+p_2(t)\lambda^2+ p_1(t) \lambda+p_0(t)$$
where $p_1(t)=q_1(t)+at$ with $a\neq 0$ and   $p_0(t), \ q_1(t),\ p_2(t)$ contains
only even powers of $t$. Then $\sige(C)\neq -\sige(C)$.
\end{lemma}

{\it Proof.} Since  rank $C_1\le 2$, $0\in \sigma(C_1)$. Therefore, $\sige(C)=\sige(C_1)\cup\sige(0_{n-k})=\sige(C_1)$.

Note that for each $t\in \IR$, $f(\lambda,t) $ is a cubic polynomial in $\lambda$ with three
non-negative real roots $\lambda_1(t)\ge \lambda_2(t)\ge \lambda_3(t)\ge 0$ and $ s_{\min}(C_1-tI_3) =\sqrt{\lambda_3(t)} $.

Without loss of generality, we may assume that $a<0$. Given
$\ve >0$, $t\in \sige(C_1)\cap \IR$ if and only if $\lambda_3(t)
<\ve^2$.
Since $\lambda_3(0)=0$ and $\underset{t\to   \infty}\lim \ \lambda_3(t)=\infty$,
there exists $t_0>0$ such that $\lambda_3( t_0)=\ve^2$.
We have $t_0\not\in \sige(C)$ and
$ f(\ve^2,t_0)=0$.
But then
$$f(\ve^2, -t_0)=f(\ve^2,t_0)-2at_0\ve^2> 0 $$
Thus,
$ \lambda_3(-t_0)<\ve^2$ implying that
$-t_0 \in \sige(C)$. So,  $t_0\in  -\sige(C)$, and thus
$\sige(C)\neq -\sige(C)$.
\qed

\begin{proposition} \label{3.2a} Let $n \ge 3$ and $A \in M_n$. The following condition are
equivalent.
\begin{itemize}
\item[{\rm (a)}] $A$ is a normal matrix with at most two distinct eigenvalues.
\item[{\rm (b)}] $\sige([A,B]) = -\sige([A,B])$ for all $B \in M_n$.
\item[{\rm (c)}] $\sige([A,B]) = -\sige([A,B])$ for all rank one nilpotent $B \in M_n$.
\end{itemize}
\end{proposition}

\it Proof. \rm
Suppose (a) holds. Then
there is a unitary $V$ and $\nu \in \IC$
such that $VAV^* - \nu I = \lambda J$ with $J =I_k \oplus -I_{n-k}$.
Then for any $B \in M_n$ such that  $VBV^* = (B_{ij})_{1 \le i, j \le 2}$
with $B_{11} \in M_k, B_{22} \in B_{22}$, we have
$$C = V[A,B]V^* = 2 \lambda
\left[\begin{matrix} O_k & B_{12} \cr - B_{21} & O_{n-k}\cr \end{matrix}\right]$$
satisfies $-C = JCJ^*$.
Thus,
$$\sige([A,B])=\sigma_{\varepsilon}([VAV^*, VBV^*])= \sigma_{\varepsilon}(-J[A,B]J^*)
= \sigma_{\varepsilon}(-[A,B]) .$$
So, condition (b) holds.

The implication (b) $\Rightarrow$ (c) is clear. To prove (c) $\Rightarrow$ (a),
we consider the contra-positive. Assume (a) is not true. We consider 2 cases.

 {\bf Case 1.} Suppose $A$
is normal with more than two distinct eigenvalues.
We may assume that
$A = \diag(a,b,c) \oplus A_2$ such that $a,\ b$ and $ c$ are distinct. If $\Re((b-a)\overline{(c-a)})\le 0$, then we have
$\Re((b-c)\overline{a-c})=\Re((b-a+a-c)\overline{a-c})=|a-c|^2-\Re((b-a)\overline{(c-a)})>0$. Thus, we may assume that  $\Re((b-a)\overline{(c-a)})>0 $ which implies that
$$|2a-(b+c)|^2=|(b-a)+(c-a)|^2>|b-a|^2+|c-a|^2>|b-c|^2\Ra \left|a-\dfrac{b+c}{2}\right| >  \dfrac{|b-c|}{2} \,.$$
Thus, by replacing $A$ with $\dfrac{2}{(b-c)} \(A-\dfrac{(b+c)}2I \)$, we may
 assume that
$A =\diag (a,1,-1) \oplus A_2$ such that $|a|>1$. Consider the rank one nilpotent
$X= {\footnotesize\begin{bmatrix}0&-\sqrt{2}&\sqrt{2}\cr
 0&-1&1\cr
  0&-1&1\cr\end{bmatrix}\oplus 0_{n-3}}$. We have
$[A,X]=C\oplus 0_{n-3}$, where
  $C= {\footnotesize\begin{bmatrix} 0& \sqrt{2} (1 - a)&
 \sqrt{2} (1 + a)\cr  0& 0& 2\cr  0& 2& 0\cr\end{bmatrix}}$.
Then
$$\det(\lambda I_3-(C-tI_3)^*(C-tI_3))=\lambda^3+p_2(t)\lambda^2+ p_1(t) \lambda+p_0(t),$$
where
$$\begin{array}{rl}
p_2(t)=&-3t^2-4|a|^2-12,\\&\\
p_1(t)=&3t^4+4\(1+|a|^2\)t^2+16\(1-|a|^2 \)t+16\(2+|a|^2 \),
\\&\\
p_0(t)=&-t^6+8t^4-16t^2 \,.\end{array}$$
Since $|a| >1$, the condition in Lemma \ref{lemt} is satisfied.
Therefore, $\sige(C) \ne -\sige(C)$.

\vskip .1in

{\bf Case 2.}
Assume that $A$ is not normal. We may assume that
$A = (a_{ij})$ is in upper triangular form such that the
$(1,2)$ entry is nonzero; see \cite[Lemma 1]{Marcus}. We may replace
$A$ by $A-a_{33}I$ and assume that $A = (A_{ij})$ with
$A_{22} \in M_{n-3}$, $A_{21} = O$, and
$$A_{11} = \begin{bmatrix}
a_{11} & a_{12} & a_{13} \cr
0 & a_{22} & a_{23} \cr
0 & 0 & 0 \cr
\end{bmatrix}.$$

{\bf Subcase (2.a)}
Suppose not both $[a_{13}, \dots, a_{1n}]$ and
$[a_{23}, \dots, a_{2n}]$ are zero.
Then there is a unitary $U = U_1 \oplus U_2$ with $U_1 \in M_2$
such that $UAU^* = \tilde A = (\tilde a_{ij})$, where
the second row of $\tilde A$ equals
$[\tilde a_{21}, \tilde a_{22}, \tilde a_{23}, 0, \dots, 0]$
with $\tilde a_{21}\in \IR$
and   $\tilde a_{21}\ne 0$ and $\tilde a_{23} \ne 0$. Let $B= E_{12}$.
Then
$$C = [\tilde A,B] = \begin{bmatrix}
-\tilde a_{21} & \tilde a_{11} - \tilde a_{22} & -\tilde a_{23} \cr
0 & \tilde a_{21} & 0 \cr
0 & 0 & 0 \cr
\end{bmatrix} \oplus O_{n-3}.$$
Then
$$\det(\lambda I_3-(C-tI_3)^*(C-tI_3))=\lambda^3+p_2(t)\lambda^2+ p_1(t) \lambda+p_0(t),$$
where
$$\begin{array}{rl}
p_2(t)=&-3t^2-|\tilde a_{22}-\tilde a_{11}|^2-|\tilde a_{23}|^2 -2\tilde a_{21}^2,\\&\\
p_1(t)=&3t^4+\(|\tilde a_{22}-\tilde a_{11}|^2 +|\tilde a_{23}|^2\)t^2-
2\tilde a_{21}|\tilde a_{23}|^2t+\tilde a_{21}^2\(\tilde a_{21}^2+|\tilde a_{23}| ^2\), \\&\\
p_0(t)=&-t^6+2\tilde a_{21}^2t^4-\tilde a_{21}^4t^2\,.\end{array}$$
Since $a_{21} $ and $\tilde a_{23} \ne 0$, the condition in Lemma \ref{lemt} is satisfied. Therefore, $\sige(C) \ne -\sige(C)$.

{\bf Subcase (2.b)} Suppose both $[a_{13}, \dots, a_{1n}]$ and
$[a_{23}, \dots, a_{2n}]$ are zero.

i) If  $a_{11}=a_{22}=0$, then we may assume that $a_{12}=1$.
Let
$$B=
\begin{bmatrix}
1 & 0 & 1 \cr
1 & 0 & 1 \cr
-1 & 0 & -1 \cr
\end{bmatrix} \oplus O_{n-3}
\ \hbox{ so that } \
C=[A,B] =
\begin{bmatrix}
1 & -1 & 1 \cr
0 & -1 & 0 \cr
0 &  1 & 0 \cr
\end{bmatrix} \oplus O_{n-3}.$$
Then
$$
\det(\lambda I_3-(C-tI_3)^*(C-tI_3))=\lambda^3+p_2(t)\lambda^2+ p_1(t) \lambda
+p_0(t),$$
where
$$\begin{array}{rl}
p_2(t)=&-3t^2-5,\\&\\
p_1(t)=&3t^4+3t^2-2t+4, \\&\\
p_0(t)=&-t^6+2t^4-t^2 \,.\end{array}$$
Therefore, the condition in Lemma \ref{lemt} is satisfied and  $\sige(C) \ne -\sige(C)$.

ii) If either $a_{11}$ or $a_{22}\neq 0$, then, applying a unitary similarity,
we may assume that $a_{11}\neq 0$. Replacing $A$ by $e^{i\theta} A$, we may assume that
$a_{11}\in \IR$. Then we may further assume that $a_{12}=1$.
Let
$B={\footnotesize \begin{bmatrix}
1 & 0 & 1 \cr
0 & 0 & 0 \cr
-1 & 0 & -1 \cr
\end{bmatrix}} \oplus O_{n-3}$ and $C=[A,B]$.
Then $C=C_1\oplus 0_{n-3} $, where
$C_1={\footnotesize\begin{bmatrix}
0 & -1 & a_{11} \cr
0 & 0 & 0 \cr
a_{11} & 1 & 0 \cr
\end{bmatrix}}.$
Then
$$
\det(\lambda I_3-(C_1-tI_3)^*(C_1-tI_3))=\lambda^3+p_2(t)\lambda^2+ p_1(t)
\lambda+p_0(t),$$
where
$$\begin{array}{rl}
p_2(t)=&-3t^2-2-2a_{11}^2,\\&\\
p_1(t)=&3t^4+2t^2-4a_{11}t+2a_{11}^2+a_{11}^4, \\&\\
p_0(t)=&-t^6+2a_{11}^2t^4-a_{11}^4t^2 \,.\end{array}$$
Therefore, the condition in Lemma \ref{lemt} is satisfied and  $\sige(C) \ne -\sige(C)$.

The proof is complete. \qed

\begin{theorem} \label{3.2} Suppose $n \ge 3$ and  $\varepsilon >0$.
Then a surjective map $\phi: M_n\rightarrow M_n$ satisfies
$$\sigma _\varepsilon([A,B])=\sigma _\varepsilon ([\Phi (A),\Phi (B)])
\qquad\hbox{ for all } A, B\in M_n$$
if and only if there exist
$\mu\in \{1, -1\}$, a unitary matrix $U\in M_n$, and
a set $\cT$ of normal matrices with at most two distinct eigenvalues
 such that
\begin{equation}\label{3.4}\phi(A) = \begin{cases} \mu U\tau(A)U^* + \nu_A I & \hbox{ if } A \in M_n \setminus \cT, \cr
-\mu U\tau(A)U^* + \nu_A I & \hbox{ if } A \in \cT, \cr \end{cases} \end{equation}
where $\nu_A \in \IC$ depends on $A$, and $\tau$ is one of the maps:
$A \mapsto A$, $A \mapsto iA^t$.
\end{theorem}

\it Proof. \rm
To prove the sufficiency, if $\tau$ has the first form,
then $
 \sige([A,B]) = \sigma_{\varepsilon}([\phi(A),\phi(B)])= \mu_A\mu_B\sige([A,B]) $
if none, one, or both of $A, B \in \cT$ by Proposition \ref{3.2a}.
If $\tau$ has the second form, then
$ \sige([A,B])= \sigma_{\varepsilon}([\phi(A),\phi(B)]) =  - \mu_A\mu_B \sige([A^t,B^t]) =
 \mu_A\mu_B\sige([A,B])$
if none, one, or both of $A, B \in \cT$ by Proposition \ref{3.2a}.

To prove the necessity,
we may compose $\phi$ by a map of the form $X \mapsto VXV^*$  and adjust
$\nu_X$ if necessary so that
$\phi$ has the form $A \mapsto \mu_A \tau(A)$, where
$\tau$ is one of the maps $A \mapsto A, A \mapsto A^t, A \mapsto \overline A,
A \mapsto A^*$. Focusing on rank one
Hermitian matrices, we see that one of the following happens.

\medskip
(1) For any rank one $A = xx^*$, $\phi(A) = \mu_A A$.
\quad
(2) For any rank one $A = xx^*$, $\phi(A) = \mu_A A^t$.

\medskip\noindent
Suppose (2) holds. We may replace $\phi$ by the map
$X \mapsto i\phi(X)^t$. Then the modified map will satisfy condition (1).
Thus, we can focus on the case when (1) holds, and prove
that $\phi$ has the asserted form with $\tau(X) = X$ for all $X \in M_n$.

In the rest of the proof, we assume that (1) holds. Then we have either

\medskip
\centerline{
i) $\phi(A) = \mu_A  A$ for all $A\in M_n$, \quad
 or \quad
ii) $\phi(A) = \mu_A  A^*$ for all $A\in M_n$.}

\medskip\noindent
We will show that i) holds with $\mu_A$ satisfying (\ref{3.4}).
Clearly, we need only consider non-scalar matrices.

\medskip
\noindent
{\bf Assertion 1} For every non-scalar matrix $A \in M_n$, $\mu_A \in \{-1,1\}$.

To prove Assertion 1,  let $A =xx^*$.
If $B = yy^*$ such that $0 \ne [A,B]$, then
$[A,B]$ is unitarily similar to $\diag(ai,-ai) \oplus O_{n-2}$
with $a = \sqrt{-\tr ([A,B]^2)/2} > 0$ so that
$$\sige([A,B]) = D(-ai,\varepsilon) \cup D(0,\varepsilon) \cup D(ai,\varepsilon).$$
Because $\sige([\phi(A),\phi(B)]) = \mu_A\mu_B \sige([A,B])$, we see that
$\mu_A\mu_B = \pm 1$.

Let   $\mu=\mu_{E_{11}}$. Suppose
$B = xx^*$ for a nonzero $x \in \IC^n$.
We can find $C = yy^*$
such that $[E_{11},C] \ne 0$ and $[B,C] \ne 0$.
Then $\mu  \mu_C, \mu_B \mu_C \in \{1, -1\}$ so that
$\mu \mu_C = \pm \mu_B \mu_C$.
It follows that $\mu_B \in \{\mu,-\mu\}$.

Choose
 $B_j = x_jx_j^*$, $j=1,\ 2$  so that $[E_{11},B_1],\ [E_{11},B_2]$ and $[B_1,B_2]\ne 0$.
Then
$$\mu \mu_{B_1},\ \mu \mu_{B_2},\  \mu_{B_2} \mu_{B_1} \in\{  1,\ -1\}.$$
Hence, $\mu^2  \in \{-1,1\}$.
So we have  either

\medskip
(a) $ \mu^2=-1\Ra \mu_B \in \{-i,i\}$  for all $B = xx^*$, \ \ or
\ \
(b) $\mu^2=1\Ra \mu_B \in \{-1,1\}$ for all $B = xx^*$.

\medskip\noindent
Next we will show that $\phi(A)=\mu_A A$   for all $A\in M_n$.
Assume the contrary that  $\phi(A)=\mu_A A^*$  for all $A\in M_n$. Let
$B_1=E_{11}+E_{13}+E_{31}+E_{33}$,
$B_2=E_{22}+E_{23}+E_{32}+E_{33}$ and  $C=E_{11}+e^{i\pi/6}E_{22}$.
Then
$$\sige([B_1,C])=D(-i,\varepsilon)\cup D(i,\varepsilon)\cup D(0,\varepsilon)$$ and
$$\sige([\phi(B_1),\phi(C)])=\mu_{B_1}\mu_{C}D(-i,\varepsilon)\cup
D(i,\varepsilon)\cup D(0,\varepsilon).$$

\medskip
\noindent
Hence, $\mu_{B_1}\mu_{C}\in\{ -1,\ 1\}.$
By a direct computation,
$$
\sige([B_2,C])=D(-e^{-2\pi i/3},\varepsilon)\cup D(e^{-2\pi i/3},\varepsilon)
\cup D(0,\varepsilon)$$
and
$$\sige([\phi(B_2),\phi(C)])=
\mu_{B_2}\mu_{C}\(D(-e^{-\pi i/3},\varepsilon)\cup D(e^{-\pi i/3},\varepsilon)
\cup D(0,\varepsilon)\).$$
Since $\mu_{B_1}=\pm \mu_{B_2}$ and  $\mu_{B_1}\mu_{C}\in\{ -1,\ 1\}$,
we have $\mu_{B_2}\mu_{C}\in\{ -1,\ 1\}$.
Hence, $\sige([\phi(B_2),\phi(C)])\neq \sige([B_2,C])$, a contradiction. Therefore, we have $\phi(A)=\mu_A A$   for all $A\in M_n$.

\medskip
For any non-scalar normal matrix $B$ with spectral
decomposition $\sum_{j=1}^n b_j x_jx_j^*$ with $b_1 \ne b_2$,
let  $C = yy^*$ with $y =x_1+ x_2 $. Then
$[B,C]$ is unitarily similar to $\diag(a,-a) \oplus O_{n-2}$.
It follows that $\mu_B\mu_C  \in \{1,-1\}$.
Because $\mu_C \mu  \in \{1,-1\}$, we see that $\mu_B \in \{\mu,-\mu\}$.
Suppose $B$ is non-normal. There is a unitary $U$ such that
$UBU^* = H + iG$, where $G = G^*$ is in diagonal form and
$H = H^*$ has a nonzero $(1,2)$ entry.
Then for $C = UE_{11}U^*$,   the matrix
$[B,C]$ is unitarily similar $\diag(a,-a) \oplus O_{n-2}$.
Again, we can conclude that $\mu_B = \pm \mu$.
So, $\mu_B \in \{\mu, -\mu\}$ for every $B \in M_n$.
Consequently, we have

\medskip\centerline{
(c) $\mu_X \in \{-i,i\}$  for all $X \in M_n$, \qquad or
\qquad
(d) $\mu_X \in \{-1,1\}$ for all $X \in M_n$.}

\bigskip\noindent
We claim that the condition (d) holds.
To this end, let $D = \diag(1,-1) \oplus O_{n-2}$ and
$B =  E_{12}/2 +  E_{23} +  E_{31}$.
Then $[D,B] =   E_{12} -  E_{23} -  E_{31}$ is a unitary matrix with
eigenvalues
$\lambda_1 =1, \lambda_2 = e^{i 2\pi/3 },
\lambda_3 = e^{i 4\pi/3 }$.
Thus,
$$\sige([D,B]) = D(\lambda_1,\varepsilon) \cup
D(\lambda_2,\varepsilon)\cup
D(\lambda_3,\varepsilon).$$
We see that $\mu_B\mu_D = 1$ for such a matrix $B$.
Similarly, if $C = -1E_{21}/2 +iE_{32} - iE_{13}$,
then $\mu_C\mu_D = 1$.
Thus, $\mu_B = \mu_C$.
Now, $[B,C] = (1+i/4) E_{11} +  (1-i/4) E_{22} -2 E_{33}$.
Then
$\mu_B\mu_C \sige([B,C]) = \sige([B,C])$ will imply that
$\mu_B \mu_C = 1$. Because, $\mu_B = \mu_C$, we see that
$\mu_B = \mu_C \in \{-1,1\}$.
Hence, the condition (d) holds.

\medskip\noindent
{\bf Assertion 2}  There is $\mu \in \{1, -1\}$ such that $\mu_A = \mu$
if  $A$ is not  a normal matrix with at most two distinct eigenvalues.

\it Proof. \rm
First we show that  for any nonzero vectors $x,\ y,\ f$  with $y,\ f\in x^\perp$, $y,\ f$
linearly independent and $\Re(f^*y)\neq 0$,
\begin{equation} \label{eq3.1}
\mu_{xf^*} = \mu_{yx^*}
\end{equation}
Note that
$C = [xf^*, y x^*]=(f^*y) x x^*
-\|x\|^2y f^*$ which has a matrix representation of the form
$$C =\left[\begin{array}{ccc}  \alpha &0 &0 \\ 0 &  -\alpha &0 \\ 0 &  \beta & 0
\end{array}\right]\oplus 0=X\oplus 0$$
with   $\alpha=f^*y\|x\|^2,\ \beta=\|x\|^2\sqrt{\|f\|^2\|y\|^2-|f^*y|^2}\neq 0$.
Then

$$\det(\lambda I_3-(X-tI_3)^*(X-tI_3))=\lambda^3+p_2(t)\lambda^2+ p_1(t) \lambda+p_0(t),$$
where
$$\begin{array}{rl}
p_2(t)=&-3t^2-(2|\alpha|^2+|\beta|^2) ,\\&\\
p_1(t)=&3t^4+(4({\rm Im}(\alpha))^2+\beta^2)t^2-
2\Re(\alpha) \beta^2 t+|\alpha|^2\( |\alpha|^2+\beta^2\), \\&\\
p_0(t)=&-t^6+(\alpha^2+\overline{\alpha}^2)t^4-|\alpha|^4t^2\,.\end{array}$$
Since $\Re(\alpha) $ and $\beta \ne 0$, the condition in Lemma \ref{lemt} is satisfied. Therefore, $\sige(C) \ne -\sige(C)$. Since
$\sige(C) = \mu_{xf^*}\mu_{yx^*}\sige(C)$,  we have  $\mu_{xf^*}\mu_{yx^*} = 1$, and thus
$\mu_{xf^*} = \mu_{yx^*}$.

If $xf^*$ and $xu^*$ are rank-1 nilpotent and if $u\in f^\perp$,
then (\ref{eq3.1})  ensures that
$$
\mu_{xf^*}=\mu_{(f+u) x^*} =\mu_{xu^*}=\mu_{u(x+f)^*}
=\mu_{f u^*} =\mu_{(x+u) f^*} =\mu_{fx^*}.
$$
So we have
\begin{equation} \label{eq3.2}
\mu_{x f^*}=\mu_{xu^*}=\mu_{f x^*}
\end{equation}
whenever $\{x,f,u\}$ is orthogonal.

Next we show that
\begin{equation}\label{3.3}
\mu_{x f^*}=\mu_{xu^*}\quad \mbox{\rm for
any nonzero vectors}\ f,u\in x^\perp.
\end{equation}

Suppose $f,\ u$ are nonzero vectors in $x^\perp$.
If $u\in f^\perp$, the equality follows from (\ref{eq3.2}).
If $u=\lambda f$ for some
nonzero scalar $\lambda$, taking $v\in \{x,f\}^\perp$ we have
$$\mu_{xf^*}=\mu_{xv^*}= \mu_{x u^*}.$$
If $u\notin f^\perp$ and $\{u,f\}$ is linearly independent, then let $v=u-cf$, where $c=\dfrac{f^*u }{f^*f}$. Then $v\in \{x,\ f\}^\perp$ and $u^*v=u^*u-\dfrac{|f^*u|^2}{f^*f}\neq 0$. By (\ref{eq3.1}) and (\ref{eq3.2}), we have
$$\mu_{xu^*}=\mu_{vx^*}=\mu_{x f^*}.$$

Next, we show that $\mu_A=\mu_B$  for
any rank one nilpotent matrices  $A,B$.
To this end,
$A=xf^*$ and
$B=yg^*$, taking unit vector $u\in\{x,y\}^\perp$ and using
({\ref{3.3}), we have
$$\mu_{xf^*}=\mu_{xu^*} = \mu_{yu^*}=\mu_{yg^*}.$$
By Proposition \ref{3.2a}, if  $A$ is not  a normal matrix with at most two distinct eigenvalues,
then there is a rank one nilpotent $B$ such that
$$-\sige([A,B]) \ne \sige([B,A]) = \mu_B\mu_A\sige([B,A]).$$
Thus, $\mu_A \mu_B  = 1$, which implies $\mu_A= \mu_B$. The desired conclusion follows.
\qed

\bigskip\noindent
{\bf Acknowledgment}

The research of the first author was partially supported by National Natural Science Foundation of
China (No.11271217).
The research of the second and third authors was partially supported by  USA NSF and
HK RCG. The second author is an honorary professor of the Shanghai University
and an honorary professor of the University of Hong Kong.

\bibliographystyle{amsplain}

\begin{thebibliography}{www}


\bibitem{CLS}
J.T. Chan, C.K. Li, and N.S. Sze,
Isometries for unitrily invariant norms, Linear Algebra Appl. 399 (2005),
53-70.

\bibitem{CLS2}
J.T. Chan, C.K. Li, and N.S. Sze,
Mappings on matrices: Invariance of functional values of matrix products,
J. of Aust. Math. Soc., 81 (2006), 165-184.

\bibitem{CFL} W.S. Cheung, S. Fallat and C.K. Li,
Multiplicative preservers on semigroups of matrices,
Liner Algebra Appl. 355 (2002), 173-186.

\bibitem{Cui1}
J. Cui, V. Forstall, C.K. Li and V. Yannello,
Properties and Preservers of the Pseudospectrum,
Linear Algebra and its Applications 436 (2012), 316-325.

\bibitem{Cui2}
J. Cui, C.K. Li and Y.T. Poon,
Pseudospectra of special operators and Pseudosectrum preservers,
J. Math. Anal. Appl. (2014), 1261-1273.

\bibitem{O1} N. Guglielmi and M.L. Overton, Fast algorithms for the approximation
of the pseudospectral abscissa and pseudospectral radius of a matrix, SIAM J. Matrix Anal. appl. 32 (2011), 1166-1192.

\bibitem{GL} R. Guralnick and C.K. Li,
Invertible preservers and algebraic groups III: preservers of unitary
similarity (congruence) invariants and overgroups of some unitary subgroups,
Linear and Multilinear Algebra 43 (1997), 257-282.

\bibitem{HJ}
R. Horn and C. Johnson, Matrix analysis.
Cambridge University Press, Cambridge, 1985.

\bibitem{Kest} H. Kestelman, Automorphisms of the field of complex numbers,
Proc. Lonon Math. Soc. 53 (1951), 1-12.

\bibitem{LPS} C.K. Li, E. Poon, and N.S. Sze,
Preservers for norms of Lie product, Operators and Matrices 3 (2009), 187-203.

\bibitem{LT} C.K. Li and N.K. Tsing,
Linear operators preserving certain classes of functions on singular values of matrices,
Linear and Multilinear Algebra 26 (1990), 133-143.

\bibitem{Marcus}
M. Marcus and M. Sandy, Conditions for the generalized numerical range to be real,
Linear Algebra Appl. 71 (1985). 219-239.

\bibitem{Sem}
P. \v Semrl,
Non-linear commutativity preserving maps.
Acta Sci. Math. (Szeged)  71  (2005), 781--819.


\end{thebibliography}

\end{document}